\documentclass[12pt, a4paper, leqno]{amsart}

\usepackage[english]{babel}

\usepackage[T1]{fontenc}
\usepackage{lmodern}

\usepackage{amsmath}
\usepackage{amssymb}
\usepackage{amsthm}

\usepackage{url}

\usepackage{hyperref} %

\usepackage{enumitem}

\newcommand{\FF}{\mathbb{F}}

\renewcommand{\epsilon}{\varepsilon}
\renewcommand{\phi}[0]{\varphi}

\newcommand{\ZZ}{\mathbb{Z}}

\newcommand{\CC}{\mathbb{C}}

\newcommand{\Order}[1]{\left\lvert #1 \right\rvert}
\newcommand{\Span}[1]{\left\langle\, #1 \,\right\rangle}

\newcommand{\Set}[1]{\left\{ #1 \right\}}

\newcommand{\norm}[0]{\trianglelefteq}

\newcommand{\thegroup}[0]{\mathfrak{A}}
\newcommand{\thealgebra}[0]{\mathfrak{g}}

\DeclareMathOperator{\ad}{ad}
\DeclareMathOperator{\Tr}{Tr}
\DeclareMathOperator{\Lie}{Lie}

\DeclareMathOperator{\ch}{ch}

\DeclareMathOperator{\GL}{GL}
\DeclareMathOperator{\SL}{SL}

\DeclareMathOperator{\SO}{SO}

\DeclareMathOperator{\Sym}{Sym}
\DeclareMathOperator{\PSL}{PSL}
\DeclareMathOperator{\PGL}{PGL}
\DeclareMathOperator{\Hom}{Hom}
\DeclareMathOperator{\Ext}{Ext}

\DeclareMathOperator{\Aut}{Aut}

\newcommand{\tce}[0]{\mathcal{G}}
\newcommand{\kt}[0]{\tilde{\kappa}}
\newcommand{\AL}{\Aut_{\mathrm{Lie}}(\thealgebra)}
\newcommand{\Hyp}[0]{\texttt{Hyp}}

\newcommand{\Gb}[0]{\mathbf{G}}

\DeclareMathOperator{\Hol}{Hol}
\DeclareMathOperator{\NHol}{NHol}

\DeclareMathOperator{\Bil}{\texttt{BilMaps}}
\DeclareMathOperator{\diag}{diag}

\DeclareMathOperator{\ind}{ind}

\newtheorem{dummy}{Dummy}
\numberwithin{dummy}{section}
\numberwithin{figure}{section}

\newtheorem{theorem}[dummy]{Theorem}

\newtheorem{lemma}[dummy]{Lemma}

\newtheorem{proposition}[dummy]{Proposition}

\newtheorem{corollary}[dummy]{Corollary} 

\theoremstyle{definition}
\newtheorem{definition}[dummy]{Definition}

\newtheorem{example}[dummy]{Example}

\theoremstyle{remark}

\newtheorem{remark}[dummy]{Remark}

\numberwithin{equation}{section}

\begin{document}

\date{26 November 2026, 11:02 CET --- Version 6.00%
} 

\title[Small multiple holomorphs]
      {Finite $p$-groups of class two with\\
      a small multiple holomorph II} 
      
\author{A. Caranti}

\address[A.~Caranti]%
 {Dipartimento di Matematica\\
  Universit\`a degli Studi di Trento\\
  via Sommarive 14\\
  I-38123 Trento\\
  Italy} 

\email{andrea.caranti@unitn.it} 

\urladdr{https://caranti.maths.unitn.it/}

\subjclass[2020]{20B35, 20D15, 20D45, 20D06, 20G05, 17B50, 16T25}

\keywords{holomorph, multiple holomorph, regular subgroups, finite
  $p$-groups of class two, equivariant bilinear maps, invariant
  trilinear forms, adjoint module, Chevalley groups, simple Lie
  algebras}

\begin{abstract}
  In a previous paper, Cindy Tsang and the author constructed finite
  $p$-groups of class two and exponent $p^{2}$, for primes $p > 2$,
  with a multiple holomorph as small as possible, in the sense that
  the quotient of the multiple holomorph by the holomorph has the
  minimal possible order $p - 1$. Those groups had small rank, and the
  methods were direct.

  In this paper we construct families of finite $p$-groups of class
  two and exponent $p^{2}$, for suitable primes $p$, of arbitrarily
  large rank and with a multiple holomorph as small as possible. The
  construction is based on the adjoint modules of simple Lie algebras
  and simple groups of Lie type, and relies on the classification of
  their invariant trilinear forms.
\end{abstract}

\thanks{The     author is a    member  of
  INdAM---GNSAGA. The author gratefully  acknowledges support from the
  Department of Mathematics of the University of Trento.}

\maketitle

\thispagestyle{empty}

\section{Introduction}

Let $\tce$ be a group, $\Sym(\tce)$ be the group of permutations on the
underlying set, and $\rho : \tce \to \Sym(\tce)$ be the right regular
representation. It is well known that the permutational holomorph
$\Hol(\tce) = N_{\Sym(\tce)}(\rho(\tce))$ is the semidirect product of
$\Aut(\tce) \le \Sym(\tce)$ by $\rho(\tce)$, and it is isomorphic to
the abstract holomorph $\Aut(\tce) \ltimes \tce$ of $\tce$.

Timothy Kohl has revived~\cite{Kohl} the study of the multiple
holomorph
\begin{equation*}
  N_{\Sym(\tce)}(N_{\Sym(\tce)}(\rho(\tce)))
  = 
  N_{\Sym(\tce)}(\Hol(\tce)).
\end{equation*}
(See also~\cite{abelian, perfect, p4}.) This group acts transitively by
conjugation on the 
set of regular 
subgroups of $\Sym(\tce)$ that are normal in the permutational
holomorph of $\tce$, and are isomorphic to $\tce$. Since the
stabiliser of $\rho(\tce)$ is $\Hol(\tce)$, the quotient
\begin{equation*}
  T(\tce) = N_{\Sym(\tce)}(\Hol(\tce)) / \Hol(\tce) 
\end{equation*}
acts regularly on this set. Regular subgroups of the permutational
holomorph of $\tce$ correspond to skew brace structures whose additive
group is $\tce$~\cite{SKEW}. Normal regular subgroups of the
permutational holomorph yield instances of bi-skew
braces~\cite{bi-Childs, bi-me, bi-CS, bi-SteTra}.

In~\cite{large}, Cindy Tsang and the author proved that $T(\tce)$ can
be arbitrarily large, in the sense that given an arbitrary finite
group $\mathcal{X}$, one can find a finite $p$-group $\tce$ of class
two, for $p$ an odd prime, such that $T(\tce)$ contains an isomorphic
copy of $\mathcal{X}$.

In~\cite{p4}, the author provided a lower bound for $T(\tce)$ when
$\tce$ is a finite $p$-group of nilpotence class two, for $p$ an odd
prime, by showing that $T(\tce)$ contains a cyclic subgroup of order
$(p - 1) p^{r - 1}$, where $p^{r} = \exp(\tce')$, and hence of order
$p - 1$ when $\exp(\tce') = p$. In the same paper it is shown that if
$\tce$ is the relatively free group in the variety of
$p$-groups of class at most two and exponent $p > 2$ on $n \ge 2$ generators,
then $T(\tce)$ has order $p - 1$.

Passing from exponent $p$ to exponent $p^{2}$ is %
more
delicate, as the $p$-th power map couples $\tce/Z(\tce)$ and $\tce'$,
which %
restricts $\Aut(\tce)$ and tends to enlarge $T(\tce)$.
In~\cite{small}, Tsang and the author exhibited examples of $p$-groups
$\tce$ of class two, with $\exp(\tce) = p^{2}$ and $\exp(\tce') = p$,
such that $T(\tce)$ is as small as possible, that is, of size $p - 1$.
Those were groups of small rank, and the calculations were explicit.

The goal of this paper is to present examples of finite $p$-groups
$\tce$ of class two and exponent $p^{2}$, for any fixed prime
$p \ge 7$ (and for $p = 5$ in small Lie ranks), of arbitrarily large rank
with $T(\tce)$ of order $p - 1$. The construction is based on the
adjoint representations of simple Lie algebras and simple groups of
Lie type, and is a broad generalisation of an example in~\cite{DH}.

In Section~\ref{sec:bilinear} we describe a method, introduced
in~\cite{p4} and used in~\cite{large}, that allows one, under suitable
hypotheses, to reduce the computation of $T(\tce)$ to the
determination of certain bilinear maps
$\tce/Z(\tce) \times \tce/Z(\tce) \to Z(\tce)$ equivariant under
automorphisms.

In Section~\ref{sec:blueprint} we construct our $p$-groups $\tce$ of
class two, starting from simple Lie algebras $\thealgebra$, and show
that the group $\Aut^{c}(\tce)$ induced by $\Aut(\tce)$ on the central
quotient of $\tce$ coincides with $\Aut_{\Lie}(\thealgebra)$. For
$\thealgebra$ of type $A_{1}$, this gives back the example
of~\cite{DH}.

In Section~\ref{sec:trilinear} we show how one can compute the
equivariant bilinear maps via certain invariant trilinear forms. We
quote from~\cite{PRV} the dimensions of the space of these forms in
characteristic zero, and prove that the same values hold in
characteristic $p > 0$, for a suitable choice of $p$.

In Section~\ref{sec:sum} we combine these ingredients to prove our
main result (Theorem~\ref{thm:main}).

The details of the reduction modulo $p$ needed in
Section~\ref{sec:trilinear} are most likely obvious for experts in
modular representation theory. For the benefit of non-experts such as
the author, the expository Appendix~\ref{sec:appendix} sketches the
relevant arguments, mainly taken
from~\cite{jantzen_2003_representations}. 

In Appendix~\ref{sec:GAP} we report the details of a
\textsl{GAP}~\cite{GAP4} computation to cover certain small cases.

\section*{Notation}

Throughout this  paper, we  shall use the  following notation.   For a
group $\tce$, we write $\tce'=  [\tce, \tce]$ for the derived subgroup
of $\tce$,  $Z(\tce)$ for the  centre of $\tce$, $\Aut(\tce)$  for the
group of automorphisms  of $\tce$, $\Aut_c(\tce)$ for  the subgroup of
$\Aut(\tce)$ consisting of the automorphisms which induce the identity
on  $\tce/Z(\tce)$, and  $\Aut^{c}(\tce) \cong  \Aut(\tce)/\Aut_{c}(\tce)$
for  the   group  of   automorphisms  of  $\tce/Z(\tce)$   induced  by
$\Aut(\tce)$.

For $x,y\in \tce$, we write $x^{y} = y^{-1} x y$ and $[x,y] = x^{-1} x^{y}$.
We will also use $[x, y]$ for commutators in a Lie algebra, but it
should be clear from context which is which.

We make use of the following identity, which holds in any group $\tce$
of nilpotence class $2$. 
\begin{equation}
  \label{eq:d-power}
  (x y)^{d}
  =
  x^{d} y^{d} [y,x]^{\binom{d}{2}}
  \text{ for all $x,y\in \tce$ and $d\in \ZZ$}.
\end{equation}

Our actions are on the right. In  particular, in a Lie algebra context
we have $(x) \ad e_{\alpha} = [x, e_{\alpha}]$.

\section{Multiple holomorph via bilinear maps}
\label{sec:bilinear}

This    section     builds    upon~\cite[Section 2]{large},    see
also~\cite[Theorem 5.2]{perfect} and \cite[Theorem 1.2]{p4}.

Let $\tce$ be a group, $\Sym(\tce)$ the group of
permutations of the underlying set $\tce$, and let
\begin{align*}
  \rho
  :\
  &
  \tce \to \Sym(\tce)
  \\&
  g \mapsto (x \mapsto x g)
\end{align*}
be the right regular representation. The \emph{(permutational)
holomorph} of $\tce$ is $\Hol(\tce) =
N_{\Sym(\tce)}(\rho(\tce))$. This is 
isomorphic to the abstract holomorph $\Aut(\tce) \ltimes
\tce$. The \emph{multiple holomorph} of $\tce$ is
\begin{equation*}
  \NHol(\tce)
  =
  N_{\Sym(\tce)}(\Hol(\tce)).
\end{equation*}

\begin{theorem}
  \label{thm:gamma}

  Let $\tce = (\tce, \cdot)$ be any finite group.

  The following data are equivalent.
  \begin{enumerate}[label=\emph{(\arabic*)}]
  \item A normal regular subgroup $N$ of $\Hol(\tce)$.
  \item
    \label{item:equivariant}
    An anti-homomorphism
    $\gamma : \tce\rightarrow \Aut(\tce)$ such that  
    \[\gamma(x^{\beta}) = \gamma(x)^{\beta}\text{ for all $x\in \tce$ and
      $\beta\in \Aut(\tce)$}.\] 
  \item  A  group operation  $\circ$  on  $\tce$ such  that
    $\Aut(\tce)\leq 
    \Aut(\tce,\circ)$ and
    \begin{equation*}
      (xy)\circ z
      =
      (x\circ z)\cdot z^{-1}\cdot (y\circ z)
      \text{ for all $x,y,z\in \tce$},
    \end{equation*}
    where $z^{-1}$ denotes the inverse of $z$ in $\tce$.
  \end{enumerate}
\end{theorem}
\begin{remark}
  \ \endgraf
  \begin{enumerate}[label=\emph{(\arabic*)}]
  \item
    In~\ref{item:equivariant} we are using exponents both for the
    action of $\beta \in \Aut(\tce)$ on $\tce$, and for conjugacy in
    $\Aut(\tce)$.
  \item 
    If the conditions of Theorem~\ref{thm:gamma} hold true, then $(\tce, \cdot,
    \circ)$ is a (right) bi-skew brace~\cite{bi-Childs, bi-me, bi-CS,
      bi-SteTra}. 
  \end{enumerate}
\end{remark}
A result of Miller~\cite{Miller} states that if two regular subgroups of
$\Sym(\tce)$  are isomorphic, then they are conjugate in
$\Sym(\tce)$. In our special case we have the following result.
\begin{proposition}
  \label{prop:theta}
  Let $\tce$ be any finite group. Let $N$ be a normal regular subgroup
  of $\Hol(\tce)$ with corresponding operation $\circ$ on $\tce$, so
  that $N = \Set{ \nu(g) : g \in \tce }$, where
  $\nu(g) : x \mapsto x \circ g$ is the right regular representation
  of $(\tce, \circ)$.

  For any $\eta \in \Sym(\tce)$ with $1^{\eta} = 1$, the following are
  equivalent:
  \begin{enumerate}[label=\emph{(\arabic*)}]
  \item $\eta : \tce \rightarrow (\tce, \circ)$ is a group isomorphism.
  \item $\eta$ conjugates $\rho(\tce)$ to $N$, namely
    $\rho(\tce)^{\eta} = N$.
  \end{enumerate}
  When these equivalent conditions hold, we automatically have
  $\eta \in \NHol(\tce)$.
\end{proposition}

\begin{proof}
  For any $x, g \in \tce$, we have $(x g)^{\eta} = x^{\rho(g) \eta}$
  and $x^{\eta} \circ g^{\eta} = x^{\eta \nu(g^{\eta})}$. Thus $\eta$
  is an isomorphism $\tce \to (\tce, \circ)$ if and only if
  $\rho(g)^{\eta} = \nu(g^{\eta})$ for all $g \in \tce$. Since
  $1^{\eta} = 1$ and both $\rho(\tce)$ and $N$ act regularly on
  $\tce$, any element of $N$ sending $1$ to $g^{\eta}$ must equal
  $\nu(g^{\eta})$, so this holds if and only if
  $\rho(\tce)^{\eta} = N$.

  Finally, if $\rho(\tce)^{\eta} = N$, then since $N \norm \Hol(\tce)$
  we have
  \begin{equation*}
    \Hol(\tce)^{\eta}
    =
    N_{\Sym(\tce)}(\rho(\tce))^{\eta}
    =
    N_{\Sym(\tce)}(N)
    =
    \Hol(\tce),
  \end{equation*}
  where the last equality follows from $\Hol(\tce) \le N_{\Sym(\tce)}(N)$
  and 
  \begin{equation*}
    \Order{N_{\Sym(\tce)}(N)}
    =
    \Order{N_{\Sym(\tce)}(\rho(\tce))}
    =
    \Order{\Hol(\tce)}.
  \end{equation*}
  Thus $\eta \in \NHol(\tce)$.
\end{proof}

The next result is a special  case of the linear technique, introduced
in~\cite[Proposition  2.2]{p4}.    For  simplicity,  we   are  writing
$\Psi(x,y)$      instead      of      $\Psi(x      Z(\tce),y
Z(\tce))$.

\begin{theorem}
  \label{thm:bilinear-maps}

  Let $p > 2$ be a prime.
  
  Let $\tce$ be a finite $p$-group of class two  such
  that $\tce' = Z(\tce)$.

  There is a one-to-one correspondence between the following:
  \begin{enumerate}[label=\emph{(\arabic*)}]
  \item
    \label{item:gamma}
    normal regular subgroups $N$ of $\Hol(\tce)$ such that
    the corresponding $\gamma(x)$ of Theorem~\ref{thm:gamma} induces
    the identity on $\tce / Z(\tce)$, for all $x 
    \in \tce$;
    \item bilinear maps $\Psi : \tce/Z(\tce)\times
      \tce/Z(\tce)\rightarrow 
      Z(\tce)$ such that 
      \begin{equation}
        \label{eqn:bilinear}
        \Psi(x^\beta,y^\beta) = \Psi(x,y)^\beta\mbox{ for all
        }x,y\in \tce\mbox{ and }\beta\in \Aut(\tce).
      \end{equation}
  \end{enumerate}
  The correspondence is given by
  \begin{equation*}
    \Psi(x,y)
    =
    x^{-1}x^{\gamma(y)}
    \text{ for all $x,y\in \tce$}.
  \end{equation*}
  The corresponding operation $\circ$ on $\tce$ is given by
  \begin{equation}\label{eqn:circ} 
    x\circ y
    =
    x^{\gamma(y)}y
    =
    x y\Psi(x,y)
    \text{ for all $x,y\in \tce$}.
  \end{equation}
\end{theorem}

We will be interested in the following situation.
\begin{lemma}
  \label{lemma:situation}
  Let $\tce$ be a finite $p$-group of class two, for a prime
  $p > 2$.

  Suppose $V = \tce/Z(\tce)$ is elementary
  abelian, so that it can be regarded 
  as a vector space over the field with $p$ elements.

  Assume    that    there    is     a    subgroup    $\thegroup$    of
  $\Aut^{c}(\tce)$ which acts irreducibly on $V$.

  Then for every $\gamma$ as in
  Theorem~\ref{thm:gamma}\ref{item:equivariant}, $\gamma(x)$ induces 
  the identity on $\tce/Z(\tce)$ for all $x \in \tce$.
\end{lemma}

\begin{proof}
  By condition~\ref{item:equivariant}
  of Theorem~\ref{thm:gamma}, $\gamma(\tce)$ is the image of the
  $p$-group $\tce$ under an anti-homomorphism which is equivariant
  under $\Aut(\tce)$. Thus $\gamma(\tce)$ is a normal
  $p$-subgroup of 
  $\Aut(\tce)$, and its image $\Gamma$ in
  $\Aut^{c}(\tce)$ is also a normal 
  $p$-subgroup. We want to show that $\Gamma = 1$.

  Since $\Gamma$ is a (possibly trivial) $p$-group acting by
  automorphisms on the $p$-group $V \ne \Set{ 0 }$, the set $U$ of
  fixed  points of 
  $\Gamma$ on 
  $V$  is a  non-trivial  subspace  of  $V$.  Since $\Gamma  \norm
  \Aut^{c}(\tce)$,   we   have  that   $U$   is   also
  invariant   under 
  $\Aut^{c}(\tce)$,  and  thus  under its subgroup
  $\thegroup$.  Since  by  assumption 
  $\thegroup$ acts  irreducibly on $V$, we  obtain that $U =  V$, that
  is, $\Gamma = 1$, as $\Aut^{c}(\tce) \le \Aut(\tce/Z(\tce))$.
\end{proof}

In what follows, let us assume that $\tce$ is a finite 
$p$-group of class two with $\tce' = Z(\tce)$, for a prime $p > 2$, so
that Theorem~\ref{thm:bilinear-maps} 
applies. Set
\begin{equation}
  \label{eq:B}
  \begin{aligned}
    \Bil
    =
    \{\,
    &\text{bilinear maps }\Psi
    :
    \tce/Z(\tce) \times
    \tce/Z(\tce)
    \rightarrow
    Z(\tce)
    \\&
    \text{satisfying~\eqref{eqn:bilinear}}
    \,\}.
  \end{aligned}
\end{equation}
This is an abelian group under pointwise multiplication
in $Z(\tce)$. Since $\tce$ has class two, the
commutator $[\, ,\,]$ on $\tce$ 
induces a bilinear map $\tce/Z(\tce) \times \tce/Z(\tce) \to \tce'$,
and we may use it to construct elements of $\Bil$ as follows. 

\begin{example}
  \label{example:power maps}
  For each $c\in \ZZ$, we have the power bilinear map
  \begin{equation}
    \label{eq:power-maps}
    \begin{aligned}
      \Psi_{[c]}
      :\ &
      \tce/Z(\tce) \times
      \tce/Z(\tce)
      \rightarrow
      Z(\tce) 
      \\&
      (x,y) \mapsto [x,y]^{c},
    \end{aligned}
  \end{equation}
  which is alternating, and clearly satisfies (\ref{eqn:bilinear}).
\end{example}

The following special case of~\cite[Proposition 3.1]{p4} provides a lower
bound for $T(\tce)$.
\begin{lemma}
  \label{example:power maps'}
  \label{lemma:power maps}
  Let $\tce$ be a finite $p$-group of class two for a prime $p > 2$,
  with $\tce' = Z(\tce)$ and $\exp(\tce') = \exp(\tce / Z(\tce)) =  p$.

  Then $T(\tce)$ contains a cyclic subgroup of order
  $p - 1$.
\end{lemma}

Before recalling the proof of Lemma~\ref{lemma:power maps}, we state the Baer
correspondence~\cite{Baer-corr}, which is the class-two case of the
Lazard correspondence and the Baker--Campbell--Hausdorff
formulas~\cite[Ch.~9 and 10]{khukhro}.
\begin{proposition}
  \label{prop:Baer}
  Let $\mathcal{B}$ be a non-empty set.
  
  There is a one-to-one correspondence between
  \begin{enumerate}[label=\emph{(\arabic*)}]
  \item
    \label{item:Baer-group}
    nilpotent group structures of nilpotence class at most two on
    $\mathcal{B}$ such that $x \mapsto x^{2}$ is a bijection, and
  \item
    \label{item:Baer-Lie}
    nilpotent Lie ring structures of nilpotence class at most two on
    $\mathcal{B}$ such that $x \mapsto 2 x$ is a bijection.
  \end{enumerate}
  In corresponding structures, group commutators and Lie brackets are
  the same.
\end{proposition}

\begin{proof}[Sketch of proof]
  If $(\mathcal{B}, \cdot)$ is as in~\ref{item:Baer-group}, then
  \begin{equation*}
    x + y
    =
    x \cdot y \cdot [x, y]^{-1/2}
  \end{equation*}
  defines an abelian group structure $(\mathcal{B}, +)$, which becomes a
  Lie ring with the group commutator as the Lie bracket.

  Conversely, if $(\mathcal{B}, +)$ is the additive group of a Lie ring
  in~\ref{item:Baer-Lie}, and $[\cdot , \cdot]$ is the Lie bracket, then
  \begin{equation*}
    x \cdot y
    =
    x + y + \frac{1}{2} \, [x, y]
  \end{equation*}
  defines a group structure as in~\ref{item:Baer-group}, whose commutator
  coincides with the Lie bracket.
\end{proof}

\begin{proof}[Proof of Lemma~\ref{lemma:power maps}]
  By assumption we have $p = \exp(\tce') = \exp(\tce / Z(\tce))$ and $\tce'
    = Z(\tce)$, so that $\tce$ has exponent dividing $p^{2}$.
    
  For each $c \in \ZZ$, consider $\Psi_{[c]} \in \Bil$
  of~\eqref{eq:power-maps}, and let
  \begin{equation}
    \label{eq:c-power}
    x \circ_{[c]} y
    =
    x y \Psi_{[c]}(x,y)
    =
    x y [x,y]^{c}
  \end{equation}
  denote the corresponding group operation on $\tce$. Clearly $c$ is
  well defined modulo $p$.

  According to Proposition~\ref{prop:Baer}, when
  $c \equiv -1/2 \pmod{p}$ in~\eqref{eq:c-power}, the group
  $(\tce, \circ_{[c]})$ is abelian, and thus not isomorphic to $\tce$.

  When $c \not\equiv -1/2 \pmod{p}$, 
  $(\tce, \circ_{[c]})$ is isomorphic to $\tce$~\cite{p4}. To see
  this, consider the power maps
  \begin{align*}
    \eta_{d} :\ &
    \tce
    \rightarrow
    \tce
    \\
    &
    x \mapsto
    x^{d},
  \end{align*}
  where $d$ runs over the integers modulo $p^{2}$ which are coprime to
  $p$. Clearly the $\eta_{d}$ are bijections. Moreover
  \begin{align*}
    (xy)^{\eta_{d}}
    &
    = (xy)^{d}
    = x^{d} y^{d} [y,x]^{\frac{d(d-1)}{2}}
    = x^{\eta_{d}} y^{\eta_{d}} [x^{\eta_{d}}, y^{\eta_{d}}]^{\frac{-1 + d^{-1}}{2}},
  \end{align*}
  where we have appealed to~\eqref{eq:d-power}, and $d^{-1}$ denotes
  the inverse modulo $p^{2}$. It follows that 
  $\eta_{d}$ is an isomorphism $\tce \to (\tce, \circ_{[c]})$, where
  $c = \frac{-1 + d^{-1}}{2}$ runs over the integers modulo $p$ not
  congruent to $-\frac{1}{2}$ modulo $p$ as $d$ runs over the integers
  modulo $p^{2}$ coprime to $p$. By Proposition~\ref{prop:theta},
  $\eta_{d} \in \NHol(\tce)$.

  The map 
  \begin{align*}
    &(\ZZ / p^{2} \ZZ)^{\times} \to T(\tce)
    \\
    &d \mapsto \eta_{d} \Hol(\tce)
  \end{align*}
  is   a   homomorphism,    because   $\eta_{d_{1}}   \eta_{d_{2}}   =
  \eta_{d_{1}d_{2}}$. Its  kernel consists  of the $d$  with $\eta_{d}
  \in \Hol(\tce)$. Since  $\eta_{d}$ fixes $1$ and the  stabiliser of $1$
  in $\Hol(\tce)$ is  $\Aut(\tce)$, this means $\eta_{d}  \in \Aut(\tce)$, that
  is,  $d \equiv  1  {\pmod p}$  by~\eqref{eq:d-power}. Therefore  the
  image is isomorphic to $(\ZZ/p \ZZ)^{\times}$, that is, it is cyclic
  of order $p-1$.
\end{proof}

\section{The construction}
\label{sec:DH}
\label{sec:blueprint}

The following example from~\cite{DH} will be the blueprint for the
construction of this paper. Consider case (II) from Section~1 of~\cite{DH}.
This is the $p$-group, for $p > 3$, given by
\begin{equation}
  \label{eq:apres}
  \begin{aligned}
    \tce
    =
    \langle
    x_{1}, x_{2}, x_{3}
    :\
    &
    x_{1}^{p} = [x_{2}, x_{3}],
    x_{2}^{p} = [x_{3}, x_{1}],
    x_{3}^{p} = [x_{1}, x_{2}],
    \\
    &
      [[x_{i}, x_{j}], x_{k}] = 1,
      \text{ for all $i, j, k$}
      \rangle.
  \end{aligned}
\end{equation}
It is proved in~\cite{DH} that $\tce$ is a $p$-group of class two,
exponent $p^{2}$ and order $p^{6}$, with $\tce' = Z(\tce)$, such that
\begin{equation*}
  \Aut^{c}(\tce)
  =
  \Aut(\tce)/\Aut_{c}(\tce)
  \cong
  \SO(3, p) \cong \PGL(2, p).
\end{equation*}
This group contains the simple group $\Omega(3, p) = \SO(3, p)'  \cong
\PSL(2, p)$, which acts irreducibly on $\tce/Z(\tce)$.

If in the relations of~\eqref{eq:apres} we replace the $p$-th
powers $x_{i}^{p}$ by $x_{i}$ and view $x_{1} , x_{2} , x_{3}$ as the
basis of a Lie algebra over the field $\FF_{p}$, for
$p > 2$, we obtain the multiplication table
\begin{equation}
  \label{eq:quaternion}
  \Span{
    x_{1}, x_{2}, x_{3}
    :
    x_{1} = [x_{2}, x_{3}],
    x_{2} = [x_{3}, x_{1}],
    x_{3} = [x_{1}, x_{2}]
    }.
\end{equation}
This is the algebra $\mathfrak{so}_{3}(\FF_{p})$, which is isomorphic
to $\mathfrak{sl}_{2}(\FF_{p})$ for $p > 2$, and
it is the adjoint module for $\PSL(2, p)$.

Our construction reverses this procedure, and is another instance of
the Baer correspondence of Proposition~\ref{prop:Baer}.

Let $\thealgebra$ be a simple Lie algebra over $\FF_{p}$, for a prime $p > 2$.

Let $M$ be a free $\ZZ / p^{2} \ZZ$-module of rank
$\dim_{\FF_{p}}(\thealgebra)$, and fix an identification $M/pM
= \thealgebra$. Multiplication by $p$ induces an isomorphism of
$\FF_{p}$-vector spaces $\thealgebra = M/pM \xrightarrow{\sim}
pM$, $x + pM \mapsto px$. Thus the Lie bracket on $\thealgebra$
induces an alternating bilinear map 
\begin{equation*}
  \langle \cdot, \cdot \rangle : M \times M \to pM, \qquad \langle x,
  y \rangle = p[\bar{x}, \bar{y}], 
\end{equation*}
where $\bar{x} = x + pM \in \thealgebra$.

Because $\langle x, y \rangle \in pM$, we have $\langle \langle x, y
\rangle, z \rangle = 0$ for all $x, y, z \in M$, so that $\langle
\cdot , \cdot \rangle$  is a Lie bracket which turns $M$ into
a Lie ring of nilpotence class two. Proposition~\ref{prop:Baer} now yields that
\begin{equation*}
  x \cdot y
  =
  x + y + \frac{1}{2}\langle x, y \rangle
\end{equation*}
defines on $M$ the structure of a group $\tce$ of nilpotence class two, in
which the commutator equals the Lie bracket.
Since $[\thealgebra, \thealgebra] = \thealgebra$ and
$Z(\thealgebra) = 0$, we immediately have $\mathcal{G}' =
Z(\mathcal{G}) = pM \cong \thealgebra$ and
$\mathcal{G}/Z(\mathcal{G}) = M/pM = \thealgebra$.

Moreover, since powers in $\mathcal{G}$ coincide with scalar multiples
in $M$, the  isomorphism $\thealgebra = M/pM  \xrightarrow{\sim} pM =
Z(\mathcal{G})$,  $\bar{x} \mapsto  px  = x^{p}$,  is  induced by  the
$p$-th    power    map    in   $\mathcal{G}$;    it    is    therefore
$\operatorname{Aut}(\mathcal{G})$-equivariant and identifies the group
commutator $\langle x,  y \rangle = p[\bar{x}, \bar{y}]$  with the Lie
bracket on $\thealgebra$.

Note that $\tce$ has size $p^{2 \dim(\thealgebra)}$, exponent $p^{2}$
and minimal number of generators equal to $\dim(\thealgebra)$. Thus, for
$\thealgebra = \mathfrak{so}_{3}(\FF_{p}) \cong \mathfrak{sl}_{2}(\FF_{p})$,
we obtain the group of~\eqref{eq:apres} from~\cite{DH}.

The next result will be crucial in the following.
\begin{lemma}
  \label{lemma:aut-c-is-aut-lie}
  \begin{equation*}
    \Aut^{c}(\tce) = \Aut_{\mathrm{Lie}}(\thealgebra).
  \end{equation*}
\end{lemma}
\begin{proof}
  By Proposition~\ref{prop:Baer}, every automorphism of $\mathcal{G}$ is a
  $\ZZ/p^{2}\ZZ$-linear automorphism of $M$.

  Every $\beta \in \GL(\thealgebra) = \GL(M/pM)$ lifts to a
  $\ZZ/p^{2}\ZZ$-linear automorphism $\tilde{\beta} \in
  \GL_{\ZZ/p^{2}\ZZ}(M)$. For any such lift
  $\tilde{\beta}$, $\ZZ/p^{2}\ZZ$-linearity gives
  \begin{equation*}
    (x \cdot y)^{\tilde{\beta}} = x^{\tilde{\beta}} +
    y^{\tilde{\beta}} + \frac{p}{2}[\bar{x}, \bar{y}]^{\beta}, \qquad
    x^{\tilde{\beta}} \cdot y^{\tilde{\beta}} = x^{\tilde{\beta}} +
    y^{\tilde{\beta}} + \frac{p}{2}[\bar{x}^{\beta}, \bar{y}^{\beta}].
  \end{equation*}
  Thus $\tilde{\beta} \in \Aut(\mathcal{G})$ if and only if
  $[\bar{x}, \bar{y}]^{\beta} = [\bar{x}^{\beta}, \bar{y}^{\beta}]$
  for all $\bar{x}, \bar{y} \in \thealgebra$, that is,
  $\beta \in \Aut_{\mathrm{Lie}}(\thealgebra)$.
\end{proof}

\section{Trilinear invariants of the adjoint module
  in positive characteristic}
\label{sec:prime}
\label{sec:trilinear}

Fix an irreducible root system $\Phi$ of rank $n$. Let
$\thealgebra_{\CC}$ be the complex simple Lie algebra of type $\Phi$,
and let $\thealgebra_{\ZZ}$ be a Chevalley lattice in it, with
Chevalley basis $\{h_{i}\} \cup \{e_{\alpha}\}_{\alpha \in \Phi}$. For
a field $K$ put
$\thealgebra_{K} = \thealgebra_{\ZZ} \otimes_{\ZZ} K$,
and write $\thealgebra = \thealgebra_{\FF_{p}}$.

Let $\Gb$ be the simply connected Chevalley group scheme of
type $\Phi$ over $\ZZ$. For a field $K$ we write $\Gb(K)$
for its group of $K$-points and $\Gb_{K}$ for its base change
to $K$, so that $\Gb(K) = \Gb_{K}(K)$; for instance, in
type $A_{n}$ we have $\Gb(K) = \SL_{n+1}(K)$. Thus
$\Gb_{\CC}$ is the simply connected complex group of type
$\Phi$, and $\Gb(\FF_{p})$ is the finite Chevalley group of
universal type. Let $k$ be an algebraic closure of $\FF_{p}$.
For $\alpha \in \Phi$ and a field $K$, let
$x_{\alpha} \colon (K, +) \to \Gb(K)$, $t \mapsto x_{\alpha}(t)$,
be the root homomorphism attached to $e_{\alpha}$; for $t \in K$ we
call $x_{\alpha}(t) \in \Gb(K)$ a root element.

The group $\Gb(K)$ acts on $\thealgebra_{K}$ through the
adjoint representation, which factors through the central isogeny
$\varphi \colon \Gb \to \Gb_{\mathrm{ad}}$ onto the
adjoint group of type $\Phi$~\cite[Chapter~5, Corollary~1]{Steinberg}.
Following Steinberg, we write $x'_{\alpha}(t) = \varphi(x_{\alpha}(t))$;
as an element of $\GL(\thealgebra_{K})$ it is
\begin{equation}
  \label{eq:root-element}
  x'_{\alpha}(t)
  = \exp(t \ad e_{\alpha})
  = \sum_{j=0}^{3} \frac{t^{j}}{j!} (\ad e_{\alpha})^{j}.
\end{equation}
(The sum stops at $j = 3$ because $(\ad e_{\alpha})^{4} = 0$, as
$\alpha$-strings have length at most $4$. By a theorem of
Chevalley~\cite{Che55} (see~\cite[Chapter 3]{Steinberg}, \cite[Section
  4.3]{carter_1972_simple} 
and~\cite[Section~25.5]{Hum72}), the 
operators $(\ad e_{\alpha})^{j}/j!$ preserve $\thealgebra_{\ZZ}$, so
by base change they make sense on $\thealgebra_{K}$ for every field
$K$.) Thus $x_{\alpha}(t)$ acts on $\thealgebra_{K}$ as
$x'_{\alpha}(t)$.

The elements $x_{\alpha}(t)$ with $\alpha \in \Phi$ and $t \in k$
generate $\Gb(k)$. Because $\Gb$ is simply connected,
$\Gb(\FF_{p})$ is likewise generated by the $x_{\alpha}(t)$ with
$t \in \FF_{p}$~\cite[Chapter~5, Corollary~3]{Steinberg}. Let
$\thegroup \le \GL(\thealgebra)$ be the subgroup generated by
the $x'_{\alpha}(t)$ with $\alpha \in \Phi$ and $t \in \FF_{p}$; thus
$\thegroup = \varphi(\Gb(\FF_{p}))$ is the image of
$\Gb(\FF_{p})$ in $\GL(\thealgebra)$. (It can be a
proper subgroup of the group of $\FF_{p}$-points of the adjoint group
scheme, as $\PSL(2, p) < \PGL(2, p)$ shows. The arguments below use
only that $\thegroup$ is generated by the $x'_{\alpha}(t)$ with
$t \in \FF_{p}$.) Each $x'_{\alpha}(t)$ is a Lie algebra automorphism,
so $\thegroup \le \AL$.

We regard $\Gb_{k}$ as an algebraic group over $k$ in the
classical sense, identifying it with its group of points
$\Gb(k)$. For a rational $\Gb_{k}$-module $M$ we write
$M^{\Gb}$ for the subspace of vectors fixed by $\Gb(k)$.
(This agrees with the invariants of the group scheme $\Gb_{k}$
in the sense of~\cite{jantzen_2003_representations}: since $k$ is
algebraically closed and $\Gb_{k}$ is smooth, $\Gb(k)$
is dense in $\Gb_{k}$.)

\subsection{Duality}\label{subsec:duality}
The normalised Killing form $\kt$ on
$\mathfrak{g}_{\mathbb{C}}$ is defined as the unique
$\Gb_{\mathbb{C}}$-invariant symmetric bilinear form such that
\begin{equation}
  \label{eq:kt}
  \kt(e_{\alpha}, e_{-\alpha}) = 1
\end{equation}
for long roots $\alpha$.
It takes integer values on the Chevalley
lattice $\mathfrak{g}_{\ZZ}$, and by condition~\eqref{eq:kt},
it does not reduce to zero modulo $p$.

Moreover, for $\alpha \in \Phi$ and $u, v \in \thealgebra_{\ZZ}$,
both sides of
\begin{equation*}
  \kt\bigl(u^{x'_{\alpha}(t)}, v^{x'_{\alpha}(t)}\bigr) = \kt(u, v)
\end{equation*}
are polynomials in $t$ with integer coefficients,
by~\eqref{eq:root-element}. They agree for all $t \in \CC$, as $\kt$
is $\Gb_{\CC}$-invariant, so the identity holds in $\ZZ[t]$.
Reducing modulo $p$ and specialising $t$ to elements of $\FF_{p}$, we
see that $\kt$ is $\thegroup$-invariant on $\thealgebra$.

Assume $p > 3$, and
$p \nmid (n+1)$ in type $A_{n}$. Then $\mathfrak{A}$ acts irreducibly
on $\mathfrak{g}$, by Subsection~\ref{subsec:irreducible}. Since
$\kt \ne 0$ on $\mathfrak{g}$, its radical is a proper
$\mathfrak{A}$-submodule of the irreducible module $\mathfrak{g}$, so
that $\kt$ is non-degenerate on $\mathfrak{g}$.

This implies in particular that the adjoint module is isomorphic to
its dual. Let $V = (\thealgebra, +)$ be the adjoint module for the
group $\thegroup$. Consider its dual $V^{\star}$. The normalised
Killing form yields an isomorphism of vector spaces
\begin{align*}
  &V \to V^{\star}
  \\
  &x \mapsto (f_{x} : y \mapsto \kt(x, y)). 
\end{align*}
This is an isomorphism of $\thegroup$-modules. In fact,
$a \in \thegroup$ acts on $f \in V^{\star}$ by
\begin{align*}
  f^{a}(y) = f(y^{a^{-1}}).
\end{align*}
We thus have
\begin{align*}
  f_{x}^{a}(y)
  =
  f_{x}(y^{a^{-1}})
  =
  \kt(x, y^{a^{-1}})
  =
  \kt((x^{a})^{a^{-1}}, y^{a^{-1}})
  =
  \kt(x^{a}, y)
  =
  f_{x^{a}}(y).
\end{align*}
Here we have used the $\thegroup$-invariance of $\kt$ established
above.

\subsection{Equivariant bilinear maps and trilinear invariants}
\label{subsec:trilinear}

With $V = (\thealgebra, +)$ as in Subsection~\ref{subsec:duality}, put
$H = \Aut_{\Lie}(\thealgebra)$, so that $\thegroup \le H$. We are interested in
the $H$-equivariant bilinear maps $V \times V \to V$, that is, the
elements of $\Hom_{H}(V \otimes V, V)$. Since $\thegroup \le H$, and
$\kt$ is an $\thegroup$-invariant non-degenerate form, so
that $V \cong V^{\star}$ as $\thegroup$-modules by
Subsection~\ref{subsec:duality}, we have
\begin{equation}
  \label{eq:bound}
  \begin{aligned}
    \Hom_{H}(V \otimes V, V)
    &
    \subseteq
    \Hom_{\thegroup}(V \otimes V, V)
    \cong \bigl( (V^{\star})^{\otimes 2} \otimes V \bigr)^{\thegroup}
    \\&
    \cong (V^{\otimes 3})^{\thegroup}
    \cong \bigl( (V^{\star})^{\otimes 3} \bigr)^{\thegroup} .
  \end{aligned}
\end{equation}
Explicitly, the isomorphism
$\Hom_{\thegroup}(V \otimes V, V) \xrightarrow{\sim}
\bigl((V^{\star})^{\otimes 3}\bigr)^{\thegroup}$ sends an
$\thegroup$-equivariant bilinear map $\Psi : V \times V \to V$ to the
$\thegroup$-invariant trilinear form
\begin{equation*}
  \Xi_{\Psi}(X, Y, Z) = \kt(\Psi(X, Y), Z).
\end{equation*}
Since $\kt$ is associative, that is, $\kt([X, Y], Z) = \kt(X, [Y, Z])$
(this holds on $\thealgebra_{\CC}$ by $\Gb_{\CC}$-invariance, hence on
$\thealgebra_{\ZZ}$ and, by reduction, on $\thealgebra$), the Lie bracket
$[\cdot, \cdot]$ corresponds under this map to the alternating
trilinear form
\begin{equation*}
  \omega(X, Y, Z) = \kt([X, Y], Z) = \kt(X, [Y, Z]).
\end{equation*}
Moreover, if $\Psi \in \Hom_{H}(V \otimes V, V)$, then $\Xi_{\Psi}$
is invariant under every automorphism in $H$ that preserves $\kt$.

Thus we are led to study the $\thegroup$-invariant trilinear forms on
$\thealgebra$. We do this in two steps. First we compute the invariants of the
corresponding algebraic group (Theorem~\ref{thm:algebraic}). Then we
compare them with the invariants of the finite group $\thegroup$
(Proposition~\ref{prop:finite}).

Put
\begin{equation}
  \label{eq:delta}
  \delta
  =
  \dim_{\CC} (\thealgebra_{\CC}^{\otimes 3})^{\Gb_{\CC}} .
\end{equation}
Equivalently, $\delta$ is the multiplicity of the adjoint module in
$\thealgebra_{\CC} \otimes \thealgebra_{\CC}$. It is proved in~\cite{PRV} that
\begin{equation}
  \label{eq:values-of-delta}
  \delta
  =
  \begin{cases}
    2 & \text{for type $A_{n}$, with $n \ge 2$},
    \\
    1 & \text{for all other types.}
  \end{cases}
\end{equation}

The following theorem transfers to positive characteristic the results
for characteristic zero quoted above. Its proof is given in
Appendix~\ref{sec:appendix}.

\begin{theorem}
  \label{thm:algebraic}
  \label{thm:modular-PY}
  Assume $p > 3$, and $p \nmid (n+1)$ in type $A_{n}$. Let
  $\Gb_{k}$ be, as above, the simply connected simple algebraic
  group of type $\Phi$ over the algebraic closure $k$ of $\FF_{p}$, and
  let $\thealgebra_{k} = \thealgebra_{\ZZ} \otimes k$ be its Lie algebra. Then
  \begin{equation*}
    \dim_{k} (\thealgebra_{k}^{\otimes 3})^{\Gb} = \delta .
  \end{equation*}
\end{theorem}

The next step requires a stronger hypothesis on $p$:
\begin{equation}
  \tag{\Hyp}
  \label{hyp:H}
  \begin{gathered}
    p \ge 7, \text{ with } p \nmid (n+1) \text{ in type } A_{n},\\
    \text{or } p = 5 \text{ in types } A_{1}, A_{2}, A_{3}, B_{2} .
  \end{gathered}
\end{equation}
Here the general polynomial degree bound of
Subsection~\ref{subsec:AtoGk} requires $p \ge 7$ (and $p \ge 11$ in
type $G_{2}$), so it gives no information for $p = 5$, nor for
$G_{2}(7)$. To cover $G_{2}(7)$ (which makes the bound $p \ge 7$
uniform across all types) and to illustrate that the result still
holds for $p = 5$ in small ranks (namely the types with Coxeter number
$h \le 4$), we verified the five cases $A_{1}(5)$, $A_{2}(5)$,
$A_{3}(5)$, $B_{2}(5)$ and $G_{2}(7)$ by computer.
Hypothesis~\eqref{hyp:H} implies the hypothesis of
Theorem~\ref{thm:algebraic}.

\begin{proposition}
  \label{prop:finite}
  Assume~\eqref{hyp:H}.
  
  Let $\thealgebra = \thealgebra_{\FF_{p}}$, and let $\thegroup \le \Aut_{\Lie}(\thealgebra)$ be
  the subgroup generated by the elements $x'_{\alpha}(t)$, for
  $\alpha \in \Phi$ and $t \in \FF_{p}$. Then
  \begin{equation*}
    \dim_{\FF_{p}} (\thealgebra^{\otimes 3})^{\thegroup} = \delta .
  \end{equation*}
\end{proposition}

\begin{proof}
  First suppose that $p \ge 7$, and that $p \ge 11$ in type $G_{2}$.
  Subsection~\ref{subsec:AtoGk} shows that
  $\dim_{\FF_{p}}(\thealgebra^{\otimes 3})^{\thegroup} =
  \dim_{k}(\thealgebra_{k}^{\otimes 3})^{\Gb}$.
  The right-hand side equals $\delta$ by Theorem~\ref{thm:algebraic}.

  The cases of~\eqref{hyp:H} that remain are $A_{1}(5)$, $A_{2}(5)$,
  $A_{3}(5)$, $B_{2}(5)$ and $G_{2}(7)$. In these five cases the
  equality has been verified directly by computer
  (Subsection~\ref{subsec:five-cases} and Appendix~\ref{sec:GAP}).
\end{proof}

\begin{corollary}
  \label{cor:bilinear}
  \label{cor:cor}
  Assume~\eqref{hyp:H}.
  
  Then $\Hom_{\Aut_{\Lie}(\thealgebra)}(V \otimes V, V)$ is $1$-dimensional,
  spanned by the Lie bracket.
\end{corollary}

\begin{proof}
  Put $H = \Aut_{\mathrm{Lie}}(\mathfrak{g})$. The Lie bracket is
  $H$-equivariant, by the definition of $H$. By~\eqref{eq:bound} and
  Proposition~\ref{prop:finite} we have
  \begin{equation*}
    \dim_{\mathbb{F}_{p}} \Hom_{H}(V \otimes V, V)
    \le \dim_{\mathbb{F}_{p}} \Hom_{\mathfrak{A}}(V \otimes V, V)
    = \delta.
  \end{equation*}
  This settles the case $\delta = 1$.

  Now let $\mathfrak{g} = \mathfrak{sl}_{n+1}(\mathbb{F}_{p})$ be of
  type $A_{n}$, with $n \ge 2$, so that $\delta = 2$. On
  $\mathfrak{sl}_{n+1}(\CC)$, the form $(X, Y) \mapsto \Tr(XY)$ is
  $\Gb_{\CC}$-invariant, symmetric, and satisfies
  $\Tr(e_{\alpha} e_{-\alpha}) = 1$ for $e_{\alpha} = E_{i,j}$ and
  $e_{-\alpha} = E_{j,i}$ ($i \ne j$); by uniqueness of the normalised
  invariant form~\eqref{eq:kt} over $\CC$, it coincides with
  $\kt$ on $\mathfrak{g}_{\ZZ}$, and hence
  $\kt(X, Y) = \Tr(XY)$ on $\mathfrak{g}$.

  The trace form $C(X, Y, Z) =  \Tr(XYZ + XZY)$ is invariant under the
  action   of   $\GL_{n+1}(\mathbb{F}_{p})$   on   $\mathfrak{g}$   by
  conjugation, and hence under $\mathfrak{A}$, as $x_{\alpha}(t)$ acts
  on $\mathfrak{sl}_{n+1}$  as conjugation by  $I + t E_{i,j}$.  It is
  symmetric, and since  $p > 3$ it corresponds by  polarisation to the
  cubic form $\Tr(X^{3})$,  because $C(X, X, X) =  2\Tr(X^{3})$. It is
  non-zero: $C(\xi, \xi, \xi) = -12 \ne  0$ for $\xi = \diag(1, 1, -2,
  0, \dots,  0)$. The  form $\omega(X, Y,  Z) =  \kt(X, [Y,
    Z])$ is alternating and non-zero, so $\omega$ and $C$ are linearly
  independent. By  Proposition~\ref{prop:finite} and~\eqref{eq:bound},
  the   space   of   $\mathfrak{A}$-invariant   trilinear   forms   on
  $\mathfrak{g}$  has dimension  $\delta =  2$,  so it  is spanned  by
  $\omega$ and $C$.

  The map $\tau : X \mapsto -X^{t}$ is a Lie automorphism of
  $\mathfrak{g}$, so $\tau \in H$. It preserves $\kt$,
  since $\Tr\bigl((-X^{t})(-Y^{t})\bigr) = \Tr\bigl((YX)^{t}\bigr) =
  \Tr(XY)$. Hence for any $\Psi \in \Hom_{H}(V \otimes V, V)$, the
  associated trilinear form $\Xi_{\Psi}$ is also $\tau$-invariant.
  Since $\tau$ preserves both $\kt$ and the Lie bracket, it
  fixes $\omega$. On the other hand, since
  $\Tr\bigl((-X^{t})^{3}\bigr) = -\Tr(X^{3})$, $\tau$ sends $C$ to
  $-C$. For $a, b \in \mathbb{F}_{p}$, the form $a\omega + bC$ is
  $\tau$-invariant only if $2b = 0$, that is, $b = 0$, as $p$ is odd.
  It follows that $\Hom_{H}(V \otimes V, V)$ corresponds to a subspace
  of $\langle \omega \rangle$, and since the Lie bracket is
  $H$-equivariant, $\Hom_{H}(V \otimes V, V)$ is spanned by the Lie
  bracket.
\end{proof}

\begin{remark}
  In type $A_{n}$ ($n \ge 2$), the bilinear counterpart of the
  symmetric trace form $C$ %
  is 
  \begin{equation*}
    \mu(X,Y) = XY + YX - \frac{2}{n+1}\,\Tr(XY)\, I ,
  \end{equation*}
  which is well defined as $p \nmid (n+1)$, satisfies
  $\kt(\mu(X, Y), Z) = C(X, Y, Z)$ and is sent to $-\mu$ by
  $\tau$. 
\end{remark}

\section{Small \texorpdfstring{$T(\tce)$}{T(G)}}
\label{sec:sum}

\begin{theorem}\label{thm:main}
  Let $\Phi$ be an irreducible root system of rank $n$, and let $p$ be
  a prime satisfying~\eqref{hyp:H}, that is, $p \ge 7$, with
  $p \nmid (n + 1)$ if $\Phi$ is of type $A_{n}$, or $p = 5$ and
  $\Phi$ is of type $A_{1}, A_{2}, A_{3}$ or $B_{2}$.

  Let $\thealgebra = \thealgebra_{\FF_{p}}$ be the Chevalley Lie
  algebra of type $\Phi$ over $\FF_{p}$, and let $\tce$ be the finite
  $p$-group of class two constructed from $\thealgebra$ in
  Section~\ref{sec:blueprint}. Then
  \begin{equation*}
    \Order{T(\tce)} = p - 1.
  \end{equation*}
\end{theorem}

\begin{proof}
Let $\thegroup \le \AL$ be the image of the finite Chevalley group
  $\mathbf{G}(\FF_{p})$, as in Section~\ref{sec:trilinear}. By
  Subsection~\ref{subsec:irreducible}, $\thegroup$ acts irreducibly on
  $\thealgebra$. Since $\thegroup$ consists of Lie algebra
  automorphisms, $[\thealgebra, \thealgebra]$ and $Z(\thealgebra)$ are
  $\thegroup$-invariant subspaces, so each is either $0$ or
  $\thealgebra$. As $\thealgebra$ is not abelian (for instance,
  $[h_{i}, e_{\alpha_{i}}] = 2 e_{\alpha_{i}} \ne 0$ for a simple root
  $\alpha_{i}$), we get $[\thealgebra, \thealgebra] = \thealgebra$ and
  $Z(\thealgebra) = 0$. %
  Now the construction
  of Section~\ref{sec:blueprint} and
  Lemma~\ref{lemma:aut-c-is-aut-lie} apply to $\thealgebra$.

  Thus $\tce$ is a finite $p$-group of class two with $\tce' = Z(\tce)$
  and $V = \tce/Z(\tce)$ elementary abelian.

  By Lemma~\ref{lemma:aut-c-is-aut-lie}, $\Aut^{c}(\tce) = \AL \ge
  \thegroup$. Since $\thegroup$ acts irreducibly on 
  $V = \thealgebra$, Theorem~\ref{thm:bilinear-maps} and
  Lemma~\ref{lemma:situation} imply that
  the normal regular
  subgroups of $\Hol(\tce)$ correspond to the elements of the space
  $\Bil$ of~\eqref{eq:B}. Since the $p$-th power map identifies
  $\tce/Z(\tce)$ and $Z(\tce)$ as $\Aut^{c}(\tce)$-modules (both equal
  to $V = \thealgebra$), the space $\Bil$ coincides with
  $\Hom_{\AL}(V \otimes V, V)$.

  By Corollary~\ref{cor:bilinear}, the space
  $\Hom_{\AL}(V \otimes V, V)$ has dimension one in all simple types,
  and it is spanned by the Lie bracket.

  The space $\Bil$ thus consists of the $p$ maps $\Psi_{[c]}$ of
  Example~\ref{example:power maps}, for $c \in \FF_{p}$. By the proof
  of Lemma~\ref{lemma:power maps}, the maps $\Psi_{[c]}$ yield $p - 1$
  groups $(\tce, \circ_{[c]})$ isomorphic to $\tce$.

  It follows that $\Order{T(\tce)} = p - 1$.
\end{proof}

\appendix

\section{Reduction modulo \texorpdfstring{$p$}{p}}
\label{sec:appendix}

This expository section contains a sketch of the definitions and the
arguments, mainly taken from~\cite[Part~II,
  Ch.~4]{jantzen_2003_representations}, that allow one to translate
the classification of invariant trilinear forms of~\cite{PRV} from
characteristic zero to positive characteristic.

\subsection{Irreducibility of the adjoint module}
\label{subsec:irreducible}

Assume $p > 3$, and $p \nmid (n + 1)$ in type $A_{n}$. Let $K$ be a
field of characteristic $p$ (we only need $K = \FF_{p}$ and $K = k$). We
show that $\thealgebra_{K}$ is irreducible as an $\thegroup$-module; in
particular, $\thealgebra$ is irreducible under $\thegroup$, and
$\thealgebra_{k}$ is irreducible under $\thegroup$ and hence a fortiori
under $\Gb(k)$.

Let $U \subseteq \thealgebra_{K}$ be an $\thegroup$-invariant
$K$-subspace, and let $\pi : \thealgebra_{K} \to \thealgebra_{K}/U$ be
the quotient map. For $u \in U$ and $\alpha \in \Phi$, applying $\pi$
to~\eqref{eq:root-element} gives the polynomial in $t$ of degree at
most $3$
\begin{equation*}
  \pi\bigl(u\, x'_{\alpha}(t)\bigr)
  =
  \pi(u) + t\,\pi([u, e_{\alpha}])
  + \frac{t^{2}}{2}\,\pi\bigl(u (\ad e_{\alpha})^{2}\bigr)
  + \frac{t^{3}}{6}\,\pi\bigl(u (\ad e_{\alpha})^{3}\bigr).
\end{equation*}
Since $x'_{\alpha}(t) \in \thegroup$ for all $t \in \FF_{p}$ and $U$ is
$\thegroup$-invariant, this polynomial vanishes for the $p$ values
$t \in \FF_{p}$. Because $p > 3$, all of its coefficients vanish in
$\thealgebra_{K}/U$; in particular, the coefficient of $t$ vanishes, so
$[u, e_{\alpha}] \in U$. Since $h_{i} = [e_{\alpha_{i}}, e_{-\alpha_{i}}]$,
the elements $e_{\alpha}$ ($\alpha \in \Phi$) generate $\thealgebra_{K}$
as a Lie algebra, so $U$ is a Lie ideal of $\thealgebra_{K}$. Under our
assumptions on $p$, the Lie algebra $\thealgebra_{K}$ is simple
by~\cite[Corollary~2.7(a)]{Hog82}, whence $U = 0$ or
$U = \thealgebra_{K}$.

\subsection{From \texorpdfstring{$\thegroup$}{the finite group} to
  \texorpdfstring{$\Gb(k)$}{G(k)}}
\label{subsec:AtoGk}

Invariance under $\thegroup$ is expressed by a system of homogeneous
linear equations with coefficients in $\FF_{p}$, so base change to $k$
does not change the dimension of the solution space:
\begin{equation*}
  \dim_{\FF_{p}}(\thealgebra^{\otimes 3})^{\thegroup}
  =
  \dim_{k}(\thealgebra_{k}^{\otimes 3})^{\thegroup}.
\end{equation*}
Each generator $x'_{\alpha}(t)$ of $\thegroup$, with $t \in \FF_{p}$, is
the image $\varphi(x_{\alpha}(t))$ of an element of $\Gb(k)$, so
$(\thealgebra_{k}^{\otimes 3})^{\Gb(k)} \subseteq
(\thealgebra_{k}^{\otimes 3})^{\thegroup}$. We claim that equality
holds when $p \ge 7$, and $p \ge 11$ in type $G_{2}$.

Let $v \in (\thealgebra_{k}^{\otimes 3})^{\thegroup}$, and let $\alpha$
be a root. By~\eqref{eq:root-element}, the matrix entries of
$x_{\alpha}(t)$ on $\thealgebra_{k}$ are polynomials in $t$ of degree
at most $2$ (at most $3$ for the short roots of $G_{2}$, as
$(\ad e_{\alpha})^{3} = 0$ in all other cases). On
$\thealgebra_{k}^{\otimes 3}$, the map $t \mapsto v x_{\alpha}(t) - v$
is therefore a polynomial map $k \to \thealgebra_{k}^{\otimes 3}$ of
degree at most $6$ (at most $9$ in type $G_{2}$) that vanishes for the
$p$ values $t \in \FF_{p}$. When $p \ge 7$ (and $p \ge 11$ in type
$G_{2}$), it vanishes for all $t \in k$. Since the root elements
$x_{\alpha}(t)$ with $t \in k$ generate $\Gb(k)$, we obtain
$v \in (\thealgebra_{k}^{\otimes 3})^{\Gb(k)}$, and thus
\begin{equation*}
  \dim_{\FF_{p}}(\thealgebra^{\otimes 3})^{\thegroup}
  =
  \dim_{k}(\thealgebra_{k}^{\otimes 3})^{\Gb(k)}
  =
  \dim_{k}(\thealgebra_{k}^{\otimes 3})^{\Gb}.
\end{equation*}

\subsection{Five special cases}
\label{subsec:five-cases}
In the five remaining cases 
\begin{equation*}
  A_{1}(5), A_{2}(5), A_{3}(5), B_{2}(5), G_{2}(7),
\end{equation*}
a computation in \textsl{GAP}~\cite{GAP4} yielded
$\dim_{\FF_{p}}(\thealgebra^{\otimes 3})^{\thegroup} = \delta$
as well. The \textsl{GAP} code and output can be found in
Appendix~\ref{sec:GAP}.

\subsection{From characteristic  \texorpdfstring{$0$}{0} to
  characteristic \texorpdfstring{$p$}{p}} 
\label{subsec:ptozero}
We now show that, with the definition of~\eqref{eq:delta} (see
also~\eqref{eq:values-of-delta}),
\begin{equation*}
  \dim_{k}(\thealgebra_{k}^{\otimes 3})^{\Gb}
  =
  \delta,
\end{equation*}
that is, the dimension of the space of invariants does not increase
when reducing modulo $p$.

In this subsection $\Gb$ stands for the algebraic group
$\Gb_{k}$. Fix a maximal torus $T$ contained in a Borel subgroup
$B$ of $\Gb$, with character group $X(T)$ and Weyl group $W$.

For a dominant weight $\lambda$, let $\Delta(\lambda)$ be the Weyl
module (this is Jantzen's $V(\lambda)$), let
$\nabla(\lambda) = H^{0}(\lambda) = \ind_{B}^{\Gb}\lambda$
be the dual Weyl module (Jantzen's convention is that $B$ corresponds
to the \emph{negative} roots), and let $L(\lambda)$ be the simple
module of highest weight $\lambda$.
\begin{definition}
  \label{def:good}
  A $\Gb$-module $M$ has a \emph{good filtration} if it has a
  filtration $0 = M_{0} \subset M_{1} \subset \dots \subset M_{m} = M$
  with $M_{i}/M_{i-1} \cong \nabla(\lambda_{i})$ for dominant weights
  $\lambda_{i}$.
\end{definition}
Let $\ZZ[X(T)]$ denote the group algebra of $X(T)$ over
$\ZZ$, with canonical basis elements denoted by $e^{\mu}$ for
$\mu \in X(T)$, satisfying $e^{\mu} e^{\nu} = e^{\mu + \nu}$ and
$e^{0} = 1$. Any finite-dimensional rational $\Gb$-module $M$
decomposes into $T$-weight spaces
\begin{equation*}
  M
  =
  \bigoplus_{\mu \in X(T)} M_{\mu}, 
\end{equation*}
and its \emph{formal character} is defined as
\begin{equation*}
  \ch(M)
  =
  \sum_{\mu \in X(T)} \dim_{k} (M_{\mu}) e^{\mu}
  \in \ZZ[X(T)]^{W}.
\end{equation*}
Because the weight multiplicities $\dim_{k} M_{\mu}$ are recorded as
integers in $\ZZ$ (rather than via traces in $k$, which only
read dimensions modulo $p$), $\ch(M)$ recovers
$\dim_{k} (M) = \sum_{\mu} \dim_{k} (M_{\mu})$, is additive on short
exact sequences (and hence on filtrations), and satisfies
$\ch(M \otimes N) =
\ch(M)\ch(N)$. Moreover, formal characters
are preserved when passing from a $\Gb_{\CC}$-module to
its reduction modulo $p$; in particular,
$\ch(\thealgebra_{k}) =
\ch(\thealgebra_{\CC})$.

We need the following facts.
\begin{enumerate}[label=(\emph{\arabic*})]
\item
  \label{item:i}
  (\cite[II.5.10]{jantzen_2003_representations})
  $\Delta(\lambda)$ and $\nabla(\lambda)$ both have formal character
  $\chi(\lambda)$, given by Weyl's character formula and independent
  of $p$. Moreover, the $\chi(\lambda)$ for $\lambda$ dominant form a
  $\ZZ$-basis of the ring of $W$-invariant formal characters.
\item
  \label{item:ii}
  (\cite[II.2.13--2.14]{jantzen_2003_representations})
  $L(\lambda)$ is the unique simple quotient of $\Delta(\lambda)$, and
  the socle of $\nabla(\lambda)$.
\item
  \label{item:iii}
  (\cite[II.4.13 and II.4.16(a)]{jantzen_2003_representations})
  $\Ext_{\Gb}^{1}(\Delta(\lambda), \nabla(\mu)) = 0$ for all dominant
  $\lambda, \mu$, and
  \begin{equation*}
    \Hom_{\Gb}(\Delta(\lambda), \nabla(\mu))
    =
    \begin{cases}
      k & \text{if $\lambda = \mu$,}\\
      0 & \text{otherwise.}
    \end{cases}
  \end{equation*}
  For any module $M$ with a good filtration, the application of
  $\Hom_{\Gb}(\Delta(\lambda), -)$ to the filtration is exact, so that
  \begin{equation*}
    \dim_{k} \Hom_{\Gb}(\Delta(\lambda), M)
    =
    (M : \nabla(\lambda)),
  \end{equation*}
  where  $(M  :   \nabla(\lambda))$  is  the  number   of  times  that
  $\nabla(\lambda)$ occurs as a quotient in a good filtration of $M$.
\item
  \label{item:iv}
  (\cite[II.4.21]{jantzen_2003_representations}, see~\cite{Wang,
    Don85, Mat90}) The tensor product of
  two modules with a good filtration has a good filtration.
\end{enumerate}

Let $\theta$ be the highest root, which is the highest weight of the
adjoint module. By Subsection~\ref{subsec:irreducible},
$\thealgebra_{k}$ is an irreducible $\Gb$-module, so
$\thealgebra_{k} \cong L(\theta)$. Since $\thealgebra_{k}$,
$\Delta(\theta)$ and $\nabla(\theta)$ all have the same formal
character as $\thealgebra_{\CC}$, they have the same dimension;
by~\ref{item:ii}, $L(\theta)$ is a quotient of $\Delta(\theta)$ and a
submodule of $\nabla(\theta)$, so equality of dimensions gives
\begin{equation*}
  \thealgebra_{k} \cong L(\theta) \cong \Delta(\theta) \cong \nabla(\theta).
\end{equation*}

Since $\thealgebra_{k} \cong \nabla(\theta)$ has a good filtration of
length $1$, two applications of~\ref{item:iv} show that
$M = \thealgebra_{k}^{\otimes 3}$ has a good filtration. Since
$\Delta(0) \cong \nabla(0) \cong k$, we have
$M^{\Gb} = \Hom_{\Gb}(\Delta(0), M)$, and~\ref{item:iii} gives
\begin{equation*}
  \dim_{k} (\thealgebra_{k}^{\otimes 3})^{\Gb}
  = (\thealgebra_{k}^{\otimes 3} : \nabla(0)).
\end{equation*}
By~\ref{item:i}, the multiplicity $(M : \nabla(\lambda))$ is the
coefficient of $\chi(\lambda)$ when $\ch(M)$ is expanded in the basis
$\{\chi(\lambda)\}$. Since
$\ch(\thealgebra_{k}^{\otimes 3}) = (\ch \thealgebra_{k})^{3} =
(\ch \thealgebra_{\CC})^{3} = \ch(\thealgebra_{\CC}^{\otimes 3})$
and, by Weyl's complete reducibility theorem in characteristic zero,
the coefficient of $\chi(0)$ in $\ch(\thealgebra_{\CC}^{\otimes 3})$
equals $\delta$, we conclude that
\begin{equation*}
  \dim_{k}(\thealgebra_{k}^{\otimes 3})^{\Gb}
  = (\thealgebra_{k}^{\otimes 3} : \nabla(0))
  = \delta .
\end{equation*}

\section{\textsl{GAP} code and output}
\label{sec:GAP}

This section records the computations we used to deal with the five
special cases of Subsection~\ref{subsec:five-cases}.

\subsection{The code}

The following \textsl{GAP} code has been provided by Google Gemini.
{\small\begin{verbatim}
# Function to compute the dimension of
# the space of trilinear invariants.
# Self-duality allows us to compute invariants in g^{\otimes 3}
# rather than in (g^*)^{\otimes 3}.
CheckInvariants := function(type, n, p)
  local F, g, b, d, R, pos, neg, gens,
        x, ad_x, exp_x, g3, id, V, nullsp;

  F := GF(p);
  g := SimpleLieAlgebra(type, n, F);
  b := Basis(g);
  d := Length(b);

  # Get Chevalley basis elements (e_alpha, f_alpha, h_i)
  R := ChevalleyBasis(g);
  pos := R[1]; # e_alpha (positive root elements)
  neg := R[2]; # f_alpha (negative root elements)

  gens := Concatenation(pos, neg);

  id := IdentityMat(d^3, F);
  V := fail;

  for x in gens do
    # Transpose so matrices act on row vectors on the right
    ad_x := TransposedMat(AdjointMatrix(b, x));
    # exp(ad_x) = I + ad_x + (1/2)*ad_x^2 + (1/6)*ad_x^3
    exp_x := IdentityMat(d, F) + ad_x
             + (1/2)*ad_x^2 + (1/6)*ad_x^3;
    g3 := KroneckerProduct(
            KroneckerProduct(exp_x, exp_x), exp_x);

    # Intersect fixed-point subspaces iteratively:
    # v * (exp_x^{\otimes 3} - I) = 0
    if V = fail then
      V := NullspaceMat(g3 - id);
    else
      nullsp := NullspaceMat(V * g3 - V);
      if Length(nullsp) = 0 then
        V := [];
        break;
      fi;
      V := nullsp * V;
    fi;
    if Length(V) = 0 then
      break;
    fi;
  od;

  Print("Type ", type, n, " at p = ", p,
        " (dim g = ", d, "): dim (g^{\\otimes 3})^A = ",
        Length(V), "\n");
end;
\end{verbatim}}

Note that in \textsl{GAP}, \texttt{AdjointMatrix(b, x)} represents
left multiplication $y \mapsto [x, y] = -[y, x]$, so its transpose
represents $-\ad x$ in our right-action convention, and \texttt{exp\_x}
is the matrix of $x'_{\alpha}(-1)$. Since
$x'_{\alpha}(-1)^{m} = x'_{\alpha}(-m)$ for $m \in \ZZ$, the elements
$x'_{\alpha}(-1)$, for $\alpha \in \Phi$, already generate $\thegroup$,
so checking $t = -1$ suffices.

\subsection{The output}

The output was
{\small\begin{verbatim}
  # Run the 5 small-prime cases:
  gap> CheckInvariants("A", 1, 5);  # Expect: 1
  Type A1 at p = 5 (dim g = 3): dim (g^{\otimes 3})^A = 1
  gap> CheckInvariants("A", 2, 5);  # Expect: 2 (omega and C)
  Type A2 at p = 5 (dim g = 8): dim (g^{\otimes 3})^A = 2
  gap> CheckInvariants("A", 3, 5);  # Expect: 2 (omega and C)
  Type A3 at p = 5 (dim g = 15): dim (g^{\otimes 3})^A = 2
  gap> CheckInvariants("B", 2, 5);  # Expect: 1
  Type B2 at p = 5 (dim g = 10): dim (g^{\otimes 3})^A = 1
  gap> CheckInvariants("G", 2, 7);  # Expect: 1
  Type G2 at p = 7 (dim g = 14): dim (g^{\otimes 3})^A = 1
\end{verbatim}}

\section*{Acknowledgements}

The author is grateful  to  Corrado de  Concini for  his  useful advice, and
for pointing  out reference~\cite{PRV}.

The author is grateful to Willem de Graaf for several useful suggestions.

\section*{Declaration of generative AI and AI-assisted technologies in
  the manuscript preparation process} 

The author is far from being an expert in modular representation
theory, and has thus resolved to use this project as a first
experience with AI-assisted research. The author has profited from
conversations with Gemini, Claude, DeepSeek, ChatGPT, and Mistral.
In particular, Gemini was instrumental in the generalisation of the
work of~\cite{DH}, and in streamlining the construction of the finite
$p$-groups of class two in Section~\ref{sec:blueprint}.

All these tools, and Claude in particular, have been useful in
revising the drafts and in navigating modular representation theory.
Appendix~\ref{sec:appendix} relies on input from these
tools.

After using these tools, the  author reviewed and edited the content
as  needed  and takes  full  responsibility  for  the content  of  the
published article.
 
\bibliographystyle{amsalpha}
\providecommand{\bysame}{\leavevmode\hbox to3em{\hrulefill}\thinspace}
\providecommand{\MR}{\relax\ifhmode\unskip\space\fi MR }
\providecommand{\MRhref}[2]{%
  \href{http://www.ams.org/mathscinet-getitem?mr=#1}{#2}
}
\providecommand{\href}[2]{#2}

\end{document}